\documentclass[12pt]{article}
\usepackage{amsfonts}

\usepackage{amsmath}

\newtheorem{theorem}{Theorem}

\newtheorem{corollary}[theorem]{Corollary}

\newtheorem{definition}[theorem]{Definition}

\newtheorem{lemma}[theorem]{Lemma}

\newtheorem{remark}[theorem]{Remark}

\newenvironment{proof}[1][Proof]{\textbf{#1.} }{\ \rule{0.5em}{0.5em}}
\numberwithin{equation}{section}

\input{tcilatex}

\begin{document}

\title{Nabla Discrete fractional Calculus and Nabla Inequalities}
\author{George A. Anastassiou \\
Department of Mathematical Sciences\\
University of Memphis\\
Memphis, TN 38152, U.S.A.\\
ganastss@memphis.edu}
\date{}
\maketitle

\begin{abstract}
Here we define a Caputo like discrete nabla fractional difference and we
produce discrete nabla fractional Taylor formulae for the first time. We
estimate their remaiders. Then we derive related discrete nabla fractional
Opial, Ostrowski, Poincar\'{e} and Sobolev type inequalities.
\end{abstract}



\section{Introduction}

Here we follow \cite{4}.

We define the rising factorial 
\begin{equation*}
t^{\overline{n}}=t\left( t+1\right) ...\left( t+n-1\right) ,\text{ \ \ \ \ }%
n\in \mathbb{N},
\end{equation*}%
and $t^{\overline{0}}=1$. In general, let $\alpha \in \mathbb{R}$, then
define $t^{\overline{\alpha }}=\frac{\Gamma \left( t+\alpha \right) }{\Gamma
\left( t\right) }$, $t\in \mathbb{R}-\{...,-2,-1,0\}$, and $0^{\overline{%
\alpha }}=0$. Note that $\nabla \left( t^{\overline{\alpha }}\right) =\alpha
t^{\overline{\alpha -1}}$, where $\nabla y\left( t\right) =y\left( t\right)
-y\left( t-1\right) $.

For $k=2,3,...$, define $\nabla ^{k}$ inductively by $\nabla ^{k}=\nabla
\nabla ^{k-1}.$ Thus $\nabla ^{k}f\left( t\right) =\sum_{m=0}^{k}\left(
-1\right) ^{m}\left( 
\begin{array}{c}
k \\ 
m%
\end{array}%
\right) f\left( t-m\right) .$

Set $\rho \left( s\right) =s-1$, we define the $n$-th order sum of $f\left(
t\right) $ by 
\begin{equation}
\nabla _{a}^{-n}f\left( t\right) =\sum_{s=a}^{t}\frac{\left( t-\rho \left(
s\right) \right) ^{\overline{n-1}}}{\left( n-1\right) !}f\left( s\right) , 
\tag{1}  \label{1}
\end{equation}%
where $t\geq a$, $n\in \mathbb{N}$.

In general we define the $\nu $-th order fractional sum of $f$ by 
\begin{equation}
\nabla _{a}^{-\nu }f\left( t\right) =\sum_{s=a}^{t}\frac{\left( t-\rho
\left( s\right) \right) ^{\overline{\nu -1}}}{\Gamma \left( \nu \right) }%
f\left( s\right) ,  \tag{2}  \label{2}
\end{equation}%
where $\nu >0$ non-integer, $t\geq a$.

We define the fractional Caputo like nabla difference for $\mu >0$, $m-1<\mu
<m$, $m=\left\lceil \mu \right\rceil $, $\left\lceil .\right\rceil $ the
ceiling of number, $m\in \mathbb{N}$, $\nu =m-\mu $, as follows 
\begin{equation}
\nabla _{a\ast }^{\mu }f\left( t\right) =\nabla _{a}^{-\nu }\left( \nabla
^{m}f\left( t\right) \right) ,\text{ \ \ \ }t\geq a.  \tag{3}  \label{3}
\end{equation}

We mention

\begin{definition}
(\cite{4}) Here $\Delta ^{m}$ is the $m$-th order forward difference
operator, $m\in \mathbb{Z}_{+}$, 
\begin{equation*}
\left( \Delta ^{m}f\right) \left( t\right) =\sum_{k=0}^{m}\left( 
\begin{array}{c}
m \\ 
k%
\end{array}%
\right) \left( -1\right) ^{m-k}f\left( t+k\right) ,\text{ \ \ \ }t\in 
\mathbb{Z}\text{.}
\end{equation*}%
Define $t^{\left( \alpha \right) }=\frac{\Gamma \left( t+1\right) }{\Gamma
\left( t+1-\alpha \right) }$, $t\in \mathbb{R}-\{...,-2,-1\},$ $\alpha >0$,
so that

\noindent $t^{\left( n\right) }=t\left( t-1\right) ...\left( t-n+1\right) $,
for $n\in \mathbb{N}$.

Note that $t^{\overline{\alpha }}=\left( t+\alpha -1\right) ^{\left( \alpha
\right) }$.

Define for $\nu >0$ the operator 
\begin{equation}
\Delta _{a}^{-\nu }f\left( t\right) =\frac{1}{\Gamma \left( \nu \right) }%
\sum_{s=a}^{t-\nu }\left( t-s-1\right) ^{\left( \nu -1\right) }f\left(
s\right) .  \tag{4}  \label{4}
\end{equation}%
We also observe that 
\begin{equation*}
\Delta ^{m}f\left( t-m\right) =\nabla ^{m}f\left( t\right) ,\text{ \ \ \ }%
\forall \text{ }m\in \mathbb{N}\text{.}
\end{equation*}
\end{definition}

We need the law of exponents.

\begin{theorem}
(\cite{4}) Let $f$ be a real valued function, and let $\mu ,\nu >0$. Then 
\begin{equation}
\nabla _{a}^{-\nu }\left( \nabla _{a}^{-\mu }f\left( t\right) \right)
=\nabla _{a}^{-\left( \mu +\nu \right) }f\left( t\right) =\nabla _{a}^{-\mu
}\left( \nabla _{a}^{-\nu }f\left( t\right) \right) ,  \tag{5}  \label{5}
\end{equation}%
for all $t\geq a$.
\end{theorem}

We also mention the discrete Taylor formula

\begin{theorem}
\label{t3}(\cite{2}) Let $f:\mathbb{Z}\rightarrow \mathbb{R}$ be a function,
and let $a\in \mathbb{Z}$. Then, for all $t\in \mathbb{Z}$ with $t\geq a+m$,
the representation holds, 
\begin{equation}
f\left( t\right) =\sum_{k=0}^{m-1}\frac{\left( t-a\right) ^{\overline{k}}}{k!%
}\nabla ^{k}f\left( a\right) +\frac{1}{\left( m-1\right) !}\sum_{\tau
=a+1}^{t}\left( t-\tau +1\right) ^{\overline{m-1}}\nabla ^{m}f\left( \tau
\right) .  \tag{6}  \label{6}
\end{equation}
\end{theorem}

\section{Main Results}

We present the following discrete backward fractional Taylor formula

\begin{theorem}
\label{t4}Let $f:\mathbb{Z}\rightarrow \mathbb{R}$ be a function, and let $%
a\in \mathbb{Z}$. Here $m-1<\mu <m$, $m=\left\lceil \mu \right\rceil $, $\mu
>0$. Then, for all $t\in \mathbb{Z}$ with $t\geq a+m$, the representation
holds, 
\begin{equation}
f\left( t\right) =\sum_{k=0}^{m-1}\frac{\left( t-a\right) ^{\overline{k}}}{k!%
}\nabla ^{k}f\left( a\right) +\frac{1}{\Gamma \left( \mu \right) }\sum_{\tau
=a+1}^{t}\left( t-\tau +1\right) ^{\overline{\mu -1}}\nabla _{\left(
a+1\right) \ast }^{\mu }f\left( \tau \right) .  \tag{7}  \label{7}
\end{equation}
\end{theorem}

\begin{proof}
We notice that 
\begin{equation*}
\nabla _{a+1}^{-\mu }\nabla _{\left( a+1\right) \ast }^{\mu }f\left(
t\right) =\nabla _{a+1}^{-\mu }\nabla _{a+1}^{-\left( m-\mu \right) }\nabla
^{m}f\left( t\right)
\end{equation*}%
\begin{equation}
\overset{\left( \text{by (\ref{5})}\right) }{=}\nabla _{a+1}^{-\left( \mu
+m-\mu \right) }\nabla ^{m}f\left( t\right) =\nabla _{a+1}^{-m}\nabla
^{m}f\left( t\right) ,  \tag{8}  \label{8}
\end{equation}%
true for $t\geq a+1$.

But 
\begin{equation}
\nabla _{a+1}^{-m}\nabla ^{m}f\left( t\right) =\frac{1}{\left( m-1\right) !}%
\sum_{\tau =a+1}^{t}\left( t-\tau +1\right) ^{\overline{m-1}}\nabla
^{m}f\left( \tau \right) ,  \tag{9}  \label{9}
\end{equation}%
and 
\begin{equation}
\nabla _{a+1}^{-\mu }\nabla _{\left( a+1\right) \ast }^{\mu }f\left(
t\right) =\frac{1}{\Gamma \left( \mu \right) }\sum_{\tau =a+1}^{t}\left(
t-\tau +1\right) ^{\overline{\mu -1}}\nabla _{\left( a+1\right) \ast }^{\mu
}f\left( \tau \right) ,  \tag{10}  \label{10}
\end{equation}%
where $t\geq a+1$.

Then we apply Theorem \ref{t3}.

The claim is proved.
\end{proof}

\begin{corollary}
(to Theorem \ref{t4}). Additionally assume that $\nabla ^{k}f\left( a\right)
=0$, for $k=0,1,...,m-1$. Then 
\begin{equation}
f\left( t\right) =\frac{1}{\Gamma \left( \mu \right) }\sum_{\tau
=a+1}^{t}\left( t-\tau +1\right) ^{\overline{\mu -1}}\nabla _{\left(
a+1\right) \ast }^{\mu }f\left( \tau \right) ,\text{ \ \ \ }\forall \text{ }%
t\geq a+m.  \tag{11}  \label{11}
\end{equation}
\end{corollary}

We need

\begin{lemma}
(\cite{4}) Let $0\leq m-1<\nu \leq m$, $m=\left\lceil \nu \right\rceil $, $%
a\in \mathbb{N}$, $f$ defined on $\mathbb{N}_{a}=\{a,a+1,...\}$. Then 
\begin{equation}
\Delta _{a}^{-\nu }f\left( t+\nu \right) =\nabla _{a}^{-\nu }f\left(
t\right) ,\text{ \ \ \ }\forall \text{ }t\in \mathbb{N}_{a}\text{.}  \tag{12}
\label{12}
\end{equation}
\end{lemma}

\begin{theorem}
(\cite{3}) Let $p\in \mathbb{N}:\nu >p$. Then 
\begin{equation}
\Delta ^{p}\left( \Delta _{a}^{-\nu }f\left( t\right) \right) =\Delta
_{a}^{-\left( \nu -p\right) }f\left( t\right) .  \tag{13}  \label{13}
\end{equation}
\end{theorem}

We give

\begin{theorem}
\label{t8}Let $p\in \mathbb{N}:\nu >p$, $a\in \mathbb{N}$. Then 
\begin{equation}
\nabla ^{p}\left( \nabla _{\alpha }^{-\nu }f\left( t\right) \right) =\nabla
_{a}^{-\left( \nu -p\right) }f\left( t\right) ,  \tag{14}  \label{14}
\end{equation}%
for $t\in \mathbb{N}_{a}$.
\end{theorem}

\begin{proof}
We notice that 
\begin{eqnarray*}
\nabla ^{p}\left( \nabla _{a}^{-\nu }f\left( t\right) \right) &=&\Delta
^{p}\left( \nabla _{a}^{-\nu }f\right) \left( t-p\right) \\
\overset{(\text{by (\ref{12})})}{=}\Delta ^{p}\left( \Delta _{a}^{-\nu
}f\left( t-p+\nu \right) \right) &=&\left( \Delta ^{p}\Delta _{a}^{-\nu
}f\right) \left( t-p+\nu \right) =:A.
\end{eqnarray*}

Also we see that 
\begin{equation*}
\nabla _{a}^{-\left( \nu -p\right) }f\left( t\right) \overset{(\text{by (\ref%
{12})})}{=}\Delta _{a}^{-\left( \nu -p\right) }f\left( t+\nu -p\right) =:B.
\end{equation*}%
But $A=B$ by (\ref{13}), proving the claim.
\end{proof}

We make

\begin{remark}
We have 
\begin{eqnarray*}
\nabla ^{p}\left( \frac{\left( t-a\right) ^{\overline{k}}}{k!}\right)
&=&\nabla ^{p}\left( \frac{\left( t+k-1-a\right) ^{\left( k\right) }}{k!}%
\right) = \\
\Delta ^{p}\left( \frac{\left( t+k-1-a-p\right) ^{\left( k\right) }}{k!}%
\right) &=&\frac{\left( t+k-1-a-p\right) ^{\left( k-p\right) }}{\left(
k-p\right) !}=\frac{\left( t-a\right) ^{\overline{k-p}}}{\left( k-p\right) !}
\end{eqnarray*}
for $k\geq p.$

That is 
\begin{equation}
\nabla ^{p}\left( \frac{\left( t-a\right) ^{\overline{k}}}{k!}\right) =\frac{%
\left( t-a\right) ^{\overline{k-p}}}{\left( k-p\right) !},\text{ \ \ for }%
k\geq p.  \tag{15}  \label{r15}
\end{equation}
\end{remark}

We have established the following discrete backward fractional extended
Taylor's formula.

\begin{theorem}
\label{t10}Let $f:\mathbb{Z}\rightarrow \mathbb{R}$ be a function, and let $%
a\in \mathbb{Z}_{+}$. Here $m-1<\mu <m$, $m=\left\lceil \mu \right\rceil $, $%
\mu >0$. Consider $p\in \mathbb{N}:\mu >p$. Then, for all $t\geq a+m$, $t\in 
\mathbb{N}$, the representation holds, 
\begin{equation*}
\nabla ^{p}f\left( t\right) =\sum_{k=p}^{m-1}\frac{\left( t-a\right) ^{%
\overline{k-p}}}{\left( k-p\right) !}\nabla ^{k}f\left( a\right) +
\end{equation*}%
\begin{equation}
\frac{1}{\Gamma \left( \mu -p\right) }\sum_{\tau =a+1}^{t}\left( t-\tau
+1\right) ^{\overline{\mu -p-1}}\nabla _{\left( a+1\right) \ast }^{\mu
}f\left( \tau \right) .  \tag{16}  \label{16}
\end{equation}
\end{theorem}

\begin{proof}
By Theorem \ref{t8} and (\ref{r15}).
\end{proof}

\textbf{Note.} When $a\in \mathbb{Z}_{+}$, and for $p=0$ put on (\ref{16})
we get (\ref{7}).

\begin{corollary}
(to Theorem \ref{t10}). Additionally assume that $\nabla ^{k}f\left(
a\right) =0$, for $k=p,...,m-1$. Then 
\begin{equation}
\nabla ^{p}f\left( t\right) =\frac{1}{\Gamma \left( \mu -p\right) }%
\sum_{\tau =a+1}^{t}\left( t-\tau +1\right) ^{\overline{\mu -p-1}}\nabla
_{\left( a+1\right) \ast }^{\mu }f\left( \tau \right) ,\text{ \ \ }\forall 
\text{ }t\geq a+m,\text{ }t\in \mathbb{N}\text{.}  \tag{17}  \label{17}
\end{equation}
\end{corollary}

\begin{remark}
(to Theorems \ref{t4}, \ref{t10}). Let $f$ be defined on $%
[a-m+1,a-m+2,...,b] $, a discrete closed interval, where $b$ is an integer.
Then (\ref{7}) and (\ref{16}) are valid only for $t\in \lbrack a+m,b]$. Here
we must assume that $a+m<b$.
\end{remark}

\begin{remark}
We would like to find 
\begin{equation*}
\sum_{\tau =a+1}^{t}\left( t-\tau +1\right) ^{\overline{\mu -1}}=\sum_{\tau
=a+1}^{t-1}\left( t-\tau +1\right) ^{\overline{\mu -1}}+\left( 1\right) ^{%
\overline{\mu -1}}=
\end{equation*}%
\begin{equation}
\sum_{\tau =a+1}^{t-1}\left( t-\tau +1\right) ^{\overline{\mu -1}}+\Gamma
\left( \mu \right) =\sum_{\tau =a+1}^{t-1}\frac{\Gamma \left( t-\tau +\mu
\right) }{\Gamma \left( t-\tau +1\right) }+\Gamma \left( \mu \right) . 
\tag{18}  \label{18}
\end{equation}%
So still to find 
\begin{equation}
A:=\sum_{\tau =a+1}^{t-1}\frac{\Gamma \left( t-\tau +\mu \right) }{\Gamma
\left( t-\tau +1\right) }.  \tag{19}  \label{19}
\end{equation}%
We will use the following formula 
\begin{equation}
\frac{\Gamma \left( x+1\right) }{\Gamma \left( x-k+1\right) }=\frac{1}{%
\left( k+1\right) }\left( \frac{\Gamma \left( x+2\right) }{\Gamma \left(
x-k+1\right) }-\frac{\Gamma \left( x+1\right) }{\Gamma \left( x-k\right) }%
\right) ,  \tag{20}  \label{20}
\end{equation}%
where $x>k$, $x,k\in \mathbb{R}:k>-1$, $x>-1$.

So for calculating $A$ we set $x:=t-\tau +\mu -1$, $k:=\mu -1$. We observe
here that $x>-1$, $k>-1$ and $x>k$. Also we see that $x+1=t-\tau +\mu $ and $%
x-k+1=t-\tau +1$. So we have 
\begin{equation}
\frac{\Gamma \left( t-\tau +\mu \right) }{\Gamma \left( t-\tau +1\right) }=%
\frac{\Gamma \left( x+1\right) }{\Gamma \left( x-k+1\right) }=\frac{1}{\mu }%
\left( \frac{\Gamma \left( t-\tau +\mu +1\right) }{\Gamma \left( t-\tau
+1\right) }-\frac{\Gamma \left( t-\tau +\mu \right) }{\Gamma \left( t-\tau
\right) }\right) ,  \tag{21}  \label{21}
\end{equation}%
for all $\tau \in \{a+1,...,t-1\}$.

Consequently we get 
\begin{eqnarray*}
A &=&\frac{1}{\mu }\left\{ \left( \frac{\Gamma \left( t-a+\mu \right) }{%
\Gamma \left( t-a\right) }-\frac{\Gamma \left( t-a-1+\mu \right) }{\Gamma
\left( t-a-1\right) }\right) \right. + \\
&&\left( \frac{\Gamma \left( t-a-1+\mu \right) }{\Gamma \left( t-a-1\right) }%
-\frac{\Gamma \left( t-a-2+\mu \right) }{\Gamma \left( t-a-2\right) }\right)
+ \\
&&\left( \frac{\Gamma \left( t-a-2+\mu \right) }{\Gamma \left( t-a-2\right) }%
-\frac{\Gamma \left( t-a-3+\mu \right) }{\Gamma \left( t-a-3\right) }\right)
+ \\
&&... \\
&&+\left. \left( \frac{\Gamma \left( \mu +2\right) }{\Gamma \left( 2\right) }%
-\frac{\Gamma \left( \mu +1\right) }{\Gamma \left( 1\right) }\right) \right\}
\end{eqnarray*}%
(telescoping sum)%
\begin{equation}
=\frac{1}{\mu }\left\{ \frac{\Gamma \left( t-a+\mu \right) }{\Gamma \left(
t-a\right) }-\Gamma \left( \mu +1\right) \right\} =\frac{\Gamma \left(
t-a+\mu \right) }{\mu \Gamma \left( t-a\right) }-\Gamma \left( \mu \right) .
\tag{22}  \label{22}
\end{equation}%
That is 
\begin{equation}
A=\frac{\Gamma \left( t-a+\mu \right) }{\mu \Gamma \left( t-a\right) }%
-\Gamma \left( \mu \right) .  \tag{23}  \label{23}
\end{equation}%
Hence we have found that 
\begin{equation}
\sum_{\tau =a+1}^{t}\left( t-\tau +1\right) ^{\overline{\mu -1}}=\frac{%
\Gamma \left( t-a+\mu \right) }{\mu \Gamma \left( t-a\right) }=\frac{\left(
t-a\right) ^{\overline{\mu }}}{\mu }.  \tag{24}  \label{24}
\end{equation}
\end{remark}

We give

\begin{corollary}
(to Theorem \ref{t4}). We get 
\begin{equation}
\left| f\left( t\right) -\sum_{k=0}^{m-1}\frac{\left( t-a\right) ^{\overline{%
k}}}{k!}\nabla ^{k}f\left( a\right) \right| \leq \frac{\left( t-a\right) ^{%
\overline{\mu }}}{\Gamma \left( \mu +1\right) }\cdot \underset{\tau \in
\{a+1,...,t\}}{\max }\left| \nabla _{\left( a+1\right) \ast }^{\mu }f\left(
\tau \right) \right| .  \tag{25}  \label{25}
\end{equation}
\end{corollary}

\begin{proof}
Use of (\ref{7}) and (\ref{24}).
\end{proof}

\begin{corollary}
(to Theorem \ref{t10}). It holds 
\begin{equation}
\left| \nabla ^{p}f\left( t\right) -\sum_{k=p}^{m-1}\frac{\left( t-a\right)
^{\overline{k-p}}}{\left( k-p\right) !}\nabla ^{k}f\left( a\right) \right|
\leq \frac{\left( t-a\right) ^{\overline{\mu -p}}}{\Gamma \left( \mu
-p+1\right) }\cdot \underset{\tau \in \{a+1,...,t\}}{\max }\left| \nabla
_{\left( a+1\right) \ast }^{\mu }f\left( \tau \right) \right| .  \tag{26}
\label{26}
\end{equation}
\end{corollary}

\begin{proof}
Use of (\ref{16}) and (\ref{24}).
\end{proof}

We present a discrete fractional Opial inequality

\begin{theorem}
\label{t16}Let $\mu >2$, $m=\left\lceil \mu \right\rceil \geq 3;$ $p\in 
\mathbb{Z}_{+}:\mu >p;$ $a\in \mathbb{Z}_{+}$. Here $f$ is a real valued
function defined on $\{a-m+1,a-m+2,...\}$. Here $t\geq a+m$, $t\in \mathbb{N}
$. Assume that $\nabla ^{k}f\left( a\right) =0$, for $k=p,...,m-1$.

Let $\gamma ,\delta >1:\frac{1}{\gamma }+\frac{1}{\delta }=1$; $C\left( \tau
\right) >0$ for $\tau =a+1,...,t;$ and $D\left( t^{\prime }\right) \geq 0$
for $t^{\prime }=a+m,...,t$. Put 
\begin{equation}
\theta \left( t,a,\mu ,p,C,\gamma \right) :=\left( \sum_{\tau =a+1}^{t}\left[
\left( t-\tau +1\right) ^{\overline{\mu -p-1}}\left( C\left( \tau \right)
\right) ^{-1}\right] ^{\gamma }\right) ^{\frac{1}{\gamma }},\text{ \ }t\geq
a+m,  \tag{27}  \label{27}
\end{equation}%
\begin{equation}
g\left( t\right) :=\sum_{\tau =a+1}^{t}\left( C\left( \tau \right) \right)
^{\delta }\left| \nabla _{\left( a+1\right) \ast }^{\mu }f\left( \tau
\right) \right| ^{\delta }\text{, \ \ }t\geq a+1;  \tag{28}  \label{28}
\end{equation}%
\begin{equation*}
G\left( t,a,m,g\right) :=2\left( g^{2}\left( t\right) -g^{2}\left(
a+m-1\right) \right) +\frac{\left( g^{2}\left( t-1\right) -g^{2}\left(
a+m-2\right) \right) }{2}+
\end{equation*}%
\begin{equation}
2\left[ g\left( t\right) g\left( t-1\right) -g\left( a+m-1\right) g\left(
a+m-2\right) \right] ,\text{ \ }t\geq a+m.  \tag{29}  \label{29}
\end{equation}%
Call also 
\begin{equation}
K\left( t\right) :=\frac{1}{\Gamma \left( \mu -p\right) }\left(
\sum_{t^{\prime }=a+m}^{t}\left[ D\left( t^{\prime }\right) \left( C\left(
t^{\prime }\right) \right) ^{-1}\theta \left( t^{\prime },a,\mu ,p,C,\gamma
\right) \right] ^{\gamma }\right) ^{\frac{1}{\gamma }},\text{\ }t\geq a+m. 
\tag{30}  \label{30}
\end{equation}%
Then 
\begin{equation}
\sum_{t^{\prime }=a+m}^{t}D\left( t^{\prime }\right) \left| \nabla
^{p}f\left( t^{\prime }\right) \right| \left| \nabla _{\left( a+1\right)
\ast }^{\mu }f\left( t^{\prime }\right) \right| \leq K\left( t\right) \left(
G\left( t,a,m,g\right) \right) ^{\frac{1}{\delta }},  \tag{31}  \label{31}
\end{equation}%
for $t\geq a+m$.
\end{theorem}

\begin{proof}
By (\ref{17}) we have 
\begin{equation*}
\left| \nabla ^{p}f\left( t\right) \right| \leq \frac{1}{\Gamma \left( \mu
-p\right) }\sum_{\tau =a+1}^{t}\left( t-\tau +1\right) ^{\overline{\mu -p-1}%
}\left| \nabla _{\left( a+1\right) \ast }^{\mu }f\left( \tau \right) \right|
\end{equation*}%
\begin{equation*}
=\frac{1}{\Gamma \left( \mu -p\right) }\sum_{\tau =a+1}^{t}\left( t-\tau
+1\right) ^{\overline{\mu -p-1}}\left( C\left( \tau \right) \right)
^{-1}C\left( \tau \right) \left| \nabla _{\left( a+1\right) \ast }^{\mu
}f\left( \tau \right) \right|
\end{equation*}%
(by discrete H\"{o}lder's inequality) 
\begin{equation*}
\leq \frac{1}{\Gamma \left( \mu -p\right) }\left( \sum_{\tau =a+1}^{t}\left[
\left( t-\tau +1\right) ^{\overline{\mu -p-1}}\left( C\left( \tau \right)
\right) ^{-1}\right] ^{\gamma }\right) ^{\frac{1}{\gamma }}\cdot
\end{equation*}%
\begin{equation*}
\left( \sum_{\tau =a+1}^{t}\left( C\left( \tau \right) \right) ^{\delta
}\left| \nabla _{\left( a+1\right) \ast }^{\mu }f\left( \tau \right) \right|
^{\delta }\right) ^{\frac{1}{\delta }}
\end{equation*}%
\begin{equation}
=\frac{\theta \left( t,a,\mu ,p,C,\gamma \right) }{\Gamma \left( \mu
-p\right) }\left( \sum_{\tau =a+1}^{t}\left( C\left( \tau \right) \right)
^{\delta }\left| \nabla _{\left( a+1\right) \ast }^{\mu }f\left( \tau
\right) \right| ^{\delta }\right) ^{\frac{1}{\delta }}\text{, \ }\forall 
\text{ }t\geq a+m\text{.}  \tag{32}  \label{32}
\end{equation}%
We have put 
\begin{equation}
g\left( t\right) =\sum_{\tau =a+1}^{t}\left( C\left( \tau \right) \right)
^{\delta }\left| \nabla _{\left( a+1\right) \ast }^{\mu }f\left( \tau
\right) \right| ^{\delta },  \tag{33}  \label{33}
\end{equation}%
which is nondecreasing in $t\geq a+1>a-m+1$.

It holds 
\begin{equation}
\nabla g\left( t\right) =\left( C\left( t\right) \right) ^{\delta }\left|
\nabla _{\left( a+1\right) \ast }^{\mu }f\left( t\right) \right| ^{\delta },%
\text{ \ \ }t\in \{a+1,...\}.  \tag{34}  \label{34}
\end{equation}%
Hence 
\begin{equation}
\left| \nabla _{\left( a+1\right) \ast }^{\mu }f\left( t\right) \right|
=\left( \nabla g\left( t\right) \right) ^{\frac{1}{\delta }}\left( C\left(
t\right) \right) ^{-1}.  \tag{35}  \label{35}
\end{equation}%
We observe for $a+m\leq t^{\prime }\leq t$ that 
\begin{equation*}
\sum_{t^{\prime }=a+m}^{t}D\left( t^{\prime }\right) \left| \nabla
^{p}f\left( t^{\prime }\right) \right| \left| \nabla _{\left( a+1\right)
\ast }^{\mu }f\left( t^{\prime }\right) \right| \leq
\end{equation*}%
\begin{equation*}
\sum_{t^{\prime }=a+m}^{t}D\left( t^{\prime }\right) \frac{\theta \left(
t^{\prime },a,\mu ,p,C,\gamma \right) }{\Gamma \left( \mu -p\right) }\left(
g\left( t^{\prime }\right) \right) ^{\frac{1}{\delta }}\left( \nabla g\left(
t^{\prime }\right) \right) ^{\frac{1}{\delta }}\left( C\left( t^{\prime
}\right) \right) ^{-1}\leq
\end{equation*}%
(by discrete H\"{o}lder's inequality) 
\begin{equation*}
\frac{1}{\Gamma \left( \mu -p\right) }\left( \sum_{t^{\prime }=a+m}^{t}\left[
D\left( t^{\prime }\right) \left( C\left( t^{\prime }\right) \right)
^{-1}\theta \left( t^{\prime },a,\mu ,p,C,\gamma \right) \right] ^{\gamma
}\right) ^{\frac{1}{\gamma }}
\end{equation*}%
\begin{equation}
\cdot \left( \sum_{t^{\prime }=a+m}^{t}g\left( t^{\prime }\right) \cdot
\nabla g\left( t^{\prime }\right) \right) ^{\frac{1}{\delta }}.  \tag{36}
\label{36}
\end{equation}%
By $m\geq 3$ notice that $a+m-2\geq a+1$.

We define the discontinuous function 
\begin{equation*}
\psi \left( x\right) =g\left( t^{\prime }\right) +\nabla g\left( t^{\prime
}\right) \left( x-t^{\prime }+1\right) \text{, \ \ for }x\in \lbrack
t^{\prime }-1,t^{\prime }]
\end{equation*}%
a closed interval of $\mathbb{R}$, and for $t^{\prime }=a+m-1,a+m,...$ .

So $\psi \left( x\right) =g\left( t^{\prime }+1\right) +\nabla g\left(
t^{\prime }+1\right) \left( x-t^{\prime }\right) $, for $x\in \lbrack
t^{\prime },t^{\prime }+1]$, and notice that $\psi \left( t^{\prime
}-\right) =2g\left( t^{\prime }\right) -g\left( t^{\prime }-1\right) $,
while $\psi \left( t^{\prime }+\right) =g\left( t^{\prime }+1\right) $; thus 
$\psi $ in general is discontinuous. Also see that $\psi ^{\prime }\left(
x\right) =\nabla g\left( t^{\prime }\right) $, for $x\in \lbrack t^{\prime
}-1,t^{\prime }]$, for $t^{\prime }=a+m-1,...$ .

Here $g\left( t\right) ,$ $\nabla g\left( t\right) \geq 0$.

We further observe that 
\begin{equation}
g\left( t^{\prime }\right) \leq \frac{g\left( t^{\prime }\right) +\left(
2g\left( t^{\prime }\right) -g\left( t^{\prime }-1\right) \right) }{2}=\frac{%
3g\left( t^{\prime }\right) -g\left( t^{\prime }-1\right) }{2},  \tag{37}
\label{37}
\end{equation}%
for $t^{\prime }=a+m-1,...$ .

The last means that 
\begin{equation}
g\left( t^{\prime }\right) \leq \int_{t^{\prime }-1}^{t^{\prime }}\psi
\left( x\right) dx\text{, \ \ for }t^{\prime }=a+m-1,...\text{ .}  \tag{38}
\label{38}
\end{equation}%
Consequently, we get 
\begin{equation*}
g\left( t^{\prime }\right) \nabla g\left( t^{\prime }\right) \leq
\int_{t^{\prime }-1}^{t^{\prime }}\psi \left( x\right) \psi ^{\prime }\left(
x\right) dx=\int_{t^{\prime }-1}^{t^{\prime }}\psi \left( x\right) d\psi
\left( x\right) =
\end{equation*}%
\begin{equation}
\left. \frac{\left( \psi \left( x\right) \right) ^{2}}{2}\right| _{t^{\prime
}-1}^{t^{\prime }}=\frac{1}{2}\left[ \left( \psi \left( t^{\prime }\right)
\right) ^{2}-\left( \psi \left( t^{\prime }-1\right) \right) ^{2}\right] . 
\tag{39}  \label{r39}
\end{equation}%
That is 
\begin{equation}
g\left( t^{\prime }\right) \nabla g\left( t^{\prime }\right) \leq \frac{1}{2}%
\left[ \left( \psi \left( t^{\prime }\right) \right) ^{2}-\left( \psi \left(
t^{\prime }-1\right) \right) ^{2}\right] \text{, \ \ for }t^{\prime
}=a+m-1,...\text{ .}  \tag{40}  \label{40}
\end{equation}%
Therefore 
\begin{equation*}
\sum_{t^{\prime }=a+m}^{t}g\left( t^{\prime }\right) \nabla g\left(
t^{\prime }\right) \leq \frac{1}{2}\sum_{t^{\prime }=a+m}^{t}\left[ \left(
\psi \left( t^{\prime }\right) \right) ^{2}-\left( \psi \left( t^{\prime
}-1\right) \right) ^{2}\right]
\end{equation*}%
\begin{equation*}
=\frac{1}{2}\left[ \left( \left( \psi \left( a+m\right) \right) ^{2}-\left(
\psi \left( a+m-1\right) \right) ^{2}\right) +\left( \left( \psi \left(
a+m+1\right) \right) ^{2}-\left( \psi \left( a+m\right) \right) ^{2}\right)
\right.
\end{equation*}%
\begin{equation*}
+\left( \left( \psi \left( a+m+2\right) \right) ^{2}-\left( \psi \left(
a+m+1\right) \right) ^{2}\right) +...+\left. \left( \left( \psi \left(
t\right) \right) ^{2}-\left( \psi \left( t-1\right) \right) ^{2}\right) 
\right]
\end{equation*}%
\begin{equation*}
=\frac{1}{2}\left[ \left( \psi \left( t\right) \right) ^{2}-\left( \psi
\left( a+m-1\right) \right) ^{2}\right]
\end{equation*}%
\begin{equation*}
=\frac{1}{2}\left[ \left( 2g\left( t\right) -g\left( t-1\right) \right)
^{2}-\left( 2g\left( a+m-1\right) -g\left( a+m-2\right) \right) ^{2}\right]
\end{equation*}%
\begin{equation*}
=2\left( g^{2}\left( t\right) -g^{2}\left( a+m-1\right) \right) +\frac{1}{2}%
\left( g^{2}\left( t-1\right) -g^{2}\left( a+m-2\right) \right)
\end{equation*}%
\begin{equation}
-2\left[ g\left( t\right) g\left( t-1\right) -g\left( a+m-1\right) g\left(
a+m-2\right) \right] .  \tag{41}  \label{41}
\end{equation}

That is 
\begin{equation*}
\sum_{t^{\prime }=a+m}^{t}g\left( t^{\prime }\right) \nabla g\left(
t^{\prime }\right) \leq 2\left( g^{2}\left( t\right) -g^{2}\left(
a+m-1\right) \right)
\end{equation*}%
\begin{equation*}
+\frac{1}{2}\left( g^{2}\left( t-1\right) -g^{2}\left( a+m-2\right) \right) -
\end{equation*}%
\begin{equation}
2\left[ g\left( t\right) g\left( t-1\right) -g\left( a+m-1\right) g\left(
a+m-2\right) \right] ,\text{ \ }\forall \text{ }t\geq a+m.  \tag{42}
\label{42}
\end{equation}%
The last proves the claim.
\end{proof}

We give

\begin{corollary}
(to Theorem \ref{t16}). Here $f$ is a real valued function defined on $%
\{-2,-1,0,...\}$, $t\geq 3$, $t\in \mathbb{N}$. Assume $f\left( 0\right)
=f\left( -1\right) =f\left( -2\right) =0$. Put 
\begin{equation}
\overline{\theta }\left( t,2.5\right) :=\left( \sum_{\tau =1}^{t}\left[
\left( t-\tau +1\right) ^{\overline{1.5}}\right] ^{2}\right) ^{\frac{1}{2}},%
\text{ \ \ }t\geq 3,  \tag{43}  \label{43}
\end{equation}%
\begin{equation}
\overline{g}\left( t\right) :=\sum_{\tau =1}^{t}\left( \nabla _{1\ast
}^{2.5}f\left( \tau \right) \right) ^{2}\text{, \ \ }t\geq 1;  \tag{44}
\label{44}
\end{equation}%
\begin{equation*}
\overline{G}\left( t,3,\overline{g}\right) :=2\left( \overline{g}^{2}\left(
t\right) -\overline{g}^{2}\left( 2\right) \right) +\frac{\left( \overline{g}%
^{2}\left( t-1\right) -\overline{g}^{2}\left( 1\right) \right) }{2}
\end{equation*}%
\begin{equation}
+2\left[ \overline{g}\left( t\right) \overline{g}\left( t-1\right) -%
\overline{g}\left( 2\right) \overline{g}\left( 1\right) \right] \text{, \ \ }%
t\geq 3\text{.}  \tag{45}  \label{45}
\end{equation}%
Call also 
\begin{equation}
\overline{K}\left( t\right) =\frac{4}{3\sqrt{\pi }}\left( \sum_{t^{\prime
}=3}^{t}\left( \overline{\theta }\left( t^{\prime },2.5\right) \right)
^{2}\right) ^{\frac{1}{2}},\text{ \ }t\geq 3\text{.}  \tag{46}  \label{46}
\end{equation}%
Then 
\begin{equation}
\sum_{t^{\prime }=3}^{t}\left| f\left( t^{\prime }\right) \right| \left|
\nabla _{1\ast }^{2.5}f\left( t^{\prime }\right) \right| \leq \overline{K}%
\left( t\right) \left( \overline{G}\left( t,3,\overline{g}\right) \right) ^{%
\frac{1}{2}},\text{ \ \ for }t\geq 3.  \tag{47}  \label{47}
\end{equation}
\end{corollary}

\textbf{Note.} Above in (\ref{45}) we have $\overline{g}\left( 1\right) =\pi
\left( f\left( 1\right) \right) ^{2}$.

Next we present a discrete fractional nabla Ostrowski type inequality.

\begin{theorem}
Let $m-1<\mu <m$, $m=\left\lceil \mu \right\rceil $, non integerer $\mu >0$; 
$p,a\in \mathbb{Z}_{+}$ with $\mu >p$. Consider $b\in \mathbb{N}$ such that $%
a+m<b$. Let $f$ be a real valued function defined on $[a-m+1,a-m+2,...,b]$.
Here $j\in \lbrack a+m,...,b]$. Assume that $\nabla ^{k}f\left( a\right) =0$%
, for $k=p+1,...,m-1$.

Then 
\begin{equation*}
\left| \frac{1}{\left( b-a-m\right) }\sum_{j=a+m+1}^{b}\nabla ^{p}f\left(
j\right) -\nabla ^{p}f\left( a\right) \right| \leq
\end{equation*}%
\begin{equation}
\frac{\left( \left( b-a\right) ^{\overline{\mu -p+1}}-m^{\overline{\mu -p+1}%
}\right) }{\Gamma \left( \mu -p+2\right) \left( b-a-m\right) }\cdot \left( 
\underset{\tau \in \{a+1,...,b\}}{\max }\left| \nabla _{\left( a+1\right)
\ast }^{\mu }f\left( \tau \right) \right| \right) .  \tag{48}  \label{48}
\end{equation}
\end{theorem}

\begin{proof}
By (\ref{16}) we have 
\begin{equation}
\nabla ^{p}f\left( j\right) -\nabla ^{p}f\left( a\right) =\frac{1}{\Gamma
\left( \mu -p\right) }\sum_{\tau =a+1}^{j}\left( j-\tau +1\right) ^{%
\overline{\mu -p-1}}\nabla _{\left( a+1\right) \ast }^{\mu }f\left( \tau
\right) ,  \tag{49}  \label{49}
\end{equation}%
for all $j\in \lbrack a+m+1,a+m+2,...,b].$

We get that 
\begin{equation*}
\frac{1}{b-\left( a+m\right) }\sum_{j=a+m+1}^{b}\nabla ^{p}f\left( j\right)
-\nabla ^{p}f\left( a\right) =
\end{equation*}%
\begin{equation*}
\frac{1}{\left( b-a-m\right) }\sum_{j=a+m+1}^{b}\left( \nabla ^{p}f\left(
j\right) -\nabla ^{p}f\left( a\right) \right) =
\end{equation*}%
\begin{equation}
\frac{1}{\Gamma \left( \mu -p\right) \left( b-a-m\right) }%
\sum_{j=a+m+1}^{b}\left( \sum_{\tau =a+1}^{j}\left( j-\tau +1\right) ^{%
\overline{\mu -p-1}}\nabla _{\left( a+1\right) \ast }^{\mu }f\left( \tau
\right) \right) .  \tag{50}  \label{50}
\end{equation}%
Therefore we obtain 
\begin{equation*}
\left| \frac{1}{\left( b-a-m\right) }\sum_{j=a+m+1}^{b}\nabla ^{p}f\left(
j\right) -\nabla ^{p}f\left( a\right) \right| \leq 
\end{equation*}%
\begin{equation*}
\frac{1}{\Gamma \left( \mu -p\right) \left( b-a-m\right) }%
\sum_{j=a+m+1}^{b}\left( \sum_{\tau =a+1}^{j}\left( j-\tau +1\right) ^{%
\overline{\mu -p-1}}\left| \nabla _{\left( a+1\right) \ast }^{\mu }f\left(
\tau \right) \right| \right) \leq 
\end{equation*}%
\begin{equation*}
\frac{1}{\Gamma \left( \mu -p\right) \left( b-a-m\right) }\left(
\sum_{j=a+m+1}^{b}\left( \sum_{\tau =a+1}^{j}\left( j-\tau +1\right) ^{%
\overline{\mu -p-1}}\right) \right) \cdot 
\end{equation*}%
\begin{equation*}
\left( \underset{\tau \in \{a+1,...,b\}}{\max }\left| \nabla _{\left(
a+1\right) \ast }^{\mu }f\left( \tau \right) \right| \right) 
\end{equation*}%
\begin{equation*}
\overset{\text{(by (\ref{24}))}}{=}\frac{1}{\Gamma \left( \mu -p+1\right)
\left( b-a-m\right) }\left( \sum_{j=a+m+1}^{b}\left( j-a\right) ^{\overline{%
\mu -p}}\right) \cdot 
\end{equation*}%
\begin{equation*}
\left( \underset{\tau \in \{a+1,...,b\}}{\max }\left| \nabla _{\left(
a+1\right) \ast }^{\mu }f\left( \tau \right) \right| \right) 
\end{equation*}%
(by Lemma 19 of \cite{1}) 
\begin{equation*}
=\frac{1}{\Gamma \left( \mu -p+2\right) \left( b-a-m\right) }\left( \left(
b-a\right) ^{\overline{\mu -p+1}}-m^{\overline{\mu -p+1}}\right) \cdot 
\end{equation*}%
\begin{equation}
\left( \underset{\tau \in \{a+1,...,b\}}{\max }\left| \nabla _{\left(
a+1\right) \ast }^{\mu }f\left( \tau \right) \right| \right) ,  \tag{51}
\label{51}
\end{equation}%
proving the claim.
\end{proof}

Next we give a discrete nabla fractional Poincar\'{e} inequality.

\begin{theorem}
Let $\mu >p$, $p\in \mathbb{Z}_{+}$, $\mu $ non-integer, $m=\left\lceil \mu
\right\rceil ;$ $a\in \mathbb{Z}_{+}$. Here $f:[a-m+1,a-m+2,...,b]%
\rightarrow \mathbb{R}$; $a+m<b$, $b\in \mathbb{N}$, and $\nabla ^{k}f\left(
a\right) =0$, $k=p,...,m-1.$

Let $\gamma ,$ $\delta >1:\frac{1}{\gamma }+\frac{1}{\delta }=1$. Then 
\begin{equation*}
\sum_{j=a+m}^{b}\left| \nabla ^{p}f\left( j\right) \right| ^{\delta }\leq 
\frac{1}{\left( \Gamma \left( \mu -p\right) \right) ^{\delta }}\left\{
\sum_{j=a+m}^{b}\left( \sum_{\tau =a+1}^{j}\left( \left( j-\tau +1\right) ^{%
\overline{\mu -p-1}}\right) ^{\gamma }\right) ^{\frac{\delta }{\gamma }%
}\right\} \cdot 
\end{equation*}%
\begin{equation}
\left( \sum_{\tau =a+1}^{b}\left| \nabla _{\left( a+1\right) \ast }^{\mu
}f\left( \tau \right) \right| ^{\delta }\right) .  \tag{52}  \label{52}
\end{equation}
\end{theorem}

\begin{proof}
We have by (\ref{17}) that 
\begin{equation}
\nabla ^{p}f\left( j\right) =\frac{1}{\Gamma \left( \mu -p\right) }%
\sum_{\tau =a+1}^{j}\left( j-\tau +1\right) ^{\overline{\mu -p-1}}\nabla
_{\left( a+1\right) \ast }^{\mu }f\left( \tau \right) ,  \tag{53}
\label{r53}
\end{equation}
\ \ $\forall $ $j\in \lbrack a+m,a+m+1,...,b].$

Let $\gamma ,$ $\delta >1$ such that $\frac{1}{\gamma }+\frac{1}{\delta }=1$.

We observe that 
\begin{equation*}
\left| \nabla ^{p}f\left( j\right) \right| \leq \frac{1}{\Gamma \left( \mu
-p\right) }\sum_{\tau =a+1}^{j}\left( j-\tau +1\right) ^{\overline{\mu -p-1}%
}\left| \nabla _{\left( a+1\right) \ast }^{\mu }f\left( \tau \right) \right| 
\end{equation*}%
(by discrete H\"{o}lder's inequality) 
\begin{equation}
\leq \frac{1}{\Gamma \left( \mu -p\right) }\left( \sum_{\tau =a+1}^{j}\left(
\left( j-\tau +1\right) ^{\overline{\mu -p-1}}\right) ^{\gamma }\right) ^{%
\frac{1}{\gamma }}\cdot \left( \sum_{\tau =a+1}^{j}\left| \nabla _{\left(
a+1\right) \ast }^{\mu }f\left( \tau \right) \right| ^{\delta }\right) ^{%
\frac{1}{\delta }}.  \tag{54}  \label{54}
\end{equation}%
I.e. it holds 
\begin{equation*}
\left| \nabla ^{p}f\left( j\right) \right| ^{\delta }\leq \frac{1}{\left(
\Gamma \left( \mu -p\right) \right) ^{\delta }}\left( \sum_{\tau
=a+1}^{j}\left( \left( j-\tau +1\right) ^{\overline{\mu -p-1}}\right)
^{\gamma }\right) ^{\frac{\delta }{\gamma }}\cdot 
\end{equation*}%
\begin{equation*}
\left( \sum_{\tau =a+1}^{j}\left| \nabla _{\left( a+1\right) \ast }^{\mu
}f\left( \tau \right) \right| ^{\delta }\right) \leq \frac{1}{\left( \Gamma
\left( \mu -p\right) \right) ^{\delta }}\cdot 
\end{equation*}%
\begin{equation}
\left( \sum_{\tau =a+1}^{j}\left( \left( j-\tau +1\right) ^{\overline{\mu
-p-1}}\right) ^{\gamma }\right) ^{\frac{\delta }{\gamma }}\cdot \left(
\sum_{\tau =a+1}^{b}\left| \nabla _{\left( a+1\right) \ast }^{\mu }f\left(
\tau \right) \right| ^{\delta }\right) ,  \tag{55}  \label{55}
\end{equation}%
$\forall $ $j\in \lbrack a+m,b]$, a discrete interval.

Applying $\sum_{j=a+m}^{b}$ on both ends of (\ref{55}) we establish (\ref{52}%
).
\end{proof}

It follows a discrete nabla Sobolev type fractional inequality.

\begin{theorem}
\label{t20}Let $\mu >p$, $p\in \mathbb{Z}_{+}$, $\mu $ non-integer, $%
m=\left\lceil \mu \right\rceil $; $a\in \mathbb{Z}_{+}$. Here $%
f:[a-m+1,...,b]\rightarrow \mathbb{R}$; $a+m<b$, $b\in \mathbb{N}$, and $%
\nabla ^{k}f\left( a\right) =0$, $k=p,...,m-1$. Let $\gamma ,$ $\delta >1:%
\frac{1}{\gamma }+\frac{1}{\delta }=1$, and $r\geq 1$. Then 
\begin{equation*}
\left( \sum_{j=a+m}^{b}\left| \nabla ^{p}f\left( j\right) \right|
^{r}\right) ^{\frac{1}{r}}\leq \frac{1}{\Gamma \left( \mu -p\right) }\left[
\sum_{j=a+m}^{b}\left( \sum_{\tau =a+1}^{j}\left( \left( j-\tau +1\right) ^{%
\overline{\mu -p-1}}\right) ^{\gamma }\right) ^{\frac{r}{\gamma }}\right] ^{%
\frac{1}{r}}\cdot 
\end{equation*}%
\begin{equation}
\left( \sum_{\tau =a+1}^{b}\left| \nabla _{\left( a+1\right) \ast }^{\mu
}f\left( \tau \right) \right| ^{\delta }\right) ^{\frac{1}{\delta }}. 
\tag{56}  \label{56}
\end{equation}
\end{theorem}

\begin{proof}
By (\ref{54}) and $r\geq 1$ we have 
\begin{equation*}
\left| \nabla ^{p}f\left( j\right) \right| ^{r}\leq \frac{1}{\left( \Gamma
\left( \mu -p\right) \right) ^{r}}\left( \sum_{\tau =a+1}^{j}\left( \left(
j-\tau +1\right) ^{\overline{\mu -p-1}}\right) ^{\gamma }\right) ^{\frac{r}{%
\gamma }}\cdot 
\end{equation*}%
\begin{equation}
\left( \sum_{\tau =a+1}^{b}\left| \nabla _{\left( a+1\right) \ast }^{\mu
}f\left( \tau \right) \right| ^{\delta }\right) ^{\frac{r}{\delta }}\text{,
\ }\forall \text{ }j\in \lbrack a+m,...,b].  \tag{57}  \label{57}
\end{equation}%
Consequently we get 
\begin{equation*}
\sum_{j=a+m}^{b}\left| \nabla ^{p}f\left( j\right) \right| ^{r}\leq \frac{1}{%
\left( \Gamma \left( \mu -p\right) \right) ^{r}}\left[ \sum_{j=a+m}^{b}%
\left( \sum_{\tau =a+1}^{j}\left( \left( j-\tau +1\right) ^{\overline{\mu
-p-1}}\right) ^{\gamma }\right) ^{\frac{r}{\gamma }}\right] \cdot 
\end{equation*}%
\begin{equation}
\left( \sum_{\tau =a+1}^{b}\left| \nabla _{\left( a+1\right) \ast }^{\mu
}f\left( \tau \right) \right| ^{\delta }\right) ^{\frac{r}{\delta }}, 
\tag{58}  \label{58}
\end{equation}%
proving the claim.
\end{proof}

We finish with the following discrete nabla fractional average Sobolev type
inequality.

\begin{theorem}
Let $0<\mu _{1}<\mu _{2}<...<\mu _{k}$ non-integers; $m_{l}=\left\lceil \mu
_{l}\right\rceil $, $l=1,...,k,$ $k\in \mathbb{N}$. Assume $\nabla ^{\tau
}f\left( a\right) =0$, for $\tau =0,1,...,m_{k}-1$, where $%
f:[a-m_{k}+1,...,b]\rightarrow \mathbb{R}$; $b\in \mathbb{N},$ $a\in \mathbb{%
Z}_{+}$. Let $r\geq 1;$ $C_{l}\left( s\right) >0$ defined on $[a+1,...,b]$, $%
l=1,...,k;$ $a+m_{k}<b$.

Call 
\begin{equation*}
B_{l}:=\sum_{\tau =a+1}^{b}C_{l}\left( \tau \right) \left( \nabla _{\left(
a+1\right) \ast }^{\mu _{l}}f\left( \tau \right) \right) ^{2},
\end{equation*}%
\begin{equation*}
\delta ^{\ast }:=\underset{1\leq l\leq k}{\max }\left\{ \frac{1}{\left(
\Gamma \left( \mu _{l}\right) \right) ^{2}}\left[ \sum_{j=a+m_{l}}^{b}\left(
\sum_{\tau =a+1}^{j}\left( \left( j-\tau +1\right) ^{\overline{\mu _{l}-1}%
}\right) ^{2}\right) ^{\frac{r}{2}}\right] ^{\frac{2}{r}}\right\} ,
\end{equation*}%
and 
\begin{equation*}
\rho ^{\ast }:=\underset{1\leq l\leq k}{\max }\left\| \frac{1}{C_{l}\left(
\tau \right) }\right\| _{\infty ,[a+1,b]}.
\end{equation*}

Then 
\begin{equation}
\left\| f\right\| _{r,[a+m_{k},b]}\leq \sqrt{\delta ^{\ast }\rho ^{\ast }}%
\left( \frac{\sum_{l=1}^{k}B_{l}}{k}\right) ^{\frac{1}{2}}.  \tag{59}
\label{r59}
\end{equation}
\end{theorem}

\begin{proof}
We see that also $\nabla ^{\tau }f\left( a\right) =0$, $\tau =0,1,...,m_{l}-1
$, $l=1,...,k-1$. So the assumptions of Theorem \ref{t20} are fulfilled for $%
f$ and fractional orders $\mu _{l}$, $l=1,...,k$. Thus by choosing $p=0$ and 
$\gamma =\delta =2$ we apply (\ref{56}), for $l=1,...,k$, to obtain 
\begin{equation*}
\left( \sum_{j=a+m_{l}}^{b}\left| f\left( j\right) \right| ^{r}\right) ^{%
\frac{1}{r}}\leq \frac{1}{\Gamma \left( \mu _{l}\right) }\left[
\sum_{j=a+m_{l}}^{b}\left( \sum_{\tau =a+1}^{j}\left( \left( j-\tau
+1\right) ^{\overline{\mu _{l}-1}}\right) ^{2}\right) ^{\frac{r}{2}}\right]
^{\frac{1}{r}}\cdot 
\end{equation*}%
\begin{equation}
\left( \sum_{\tau =a+1}^{b}\left( \nabla _{\left( a+1\right) \ast }^{\mu
_{l}}f\left( \tau \right) \right) ^{2}\right) ^{\frac{1}{2}}.  \tag{60}
\label{60}
\end{equation}%
Hence it holds 
\begin{equation*}
\left( \sum_{j=a+m_{l}}^{b}\left| f\left( j\right) \right| ^{r}\right) ^{%
\frac{2}{r}}\leq \frac{1}{\left( \Gamma \left( \mu _{l}\right) \right) ^{2}}%
\left[ \sum_{j=a+m_{l}}^{b}\left( \sum_{\tau =a+1}^{j}\left( \left( j-\tau
+1\right) ^{\overline{\mu _{l}-1}}\right) ^{2}\right) ^{\frac{r}{2}}\right]
^{\frac{2}{r}}\cdot 
\end{equation*}%
\begin{equation*}
\left( \sum_{\tau =a+1}^{b}\left( \nabla _{\left( a+1\right) \ast }^{\mu
_{l}}f\left( \tau \right) \right) ^{2}\right) \leq \delta ^{\ast }\left(
\sum_{\tau =a+1}^{b}\left( \nabla _{\left( a+1\right) \ast }^{\mu
_{l}}f\left( \tau \right) \right) ^{2}\right) =
\end{equation*}%
\begin{equation*}
\delta ^{\ast }\left( \sum_{\tau =a+1}^{b}\left( C_{l}\left( \tau \right)
\right) ^{-1}\left( C_{l}\left( \tau \right) \right) \left( \nabla _{\left(
a+1\right) \ast }^{\mu _{l}}f\left( \tau \right) \right) ^{2}\right) \leq 
\end{equation*}%
\begin{equation*}
\delta ^{\ast }\rho ^{\ast }\left( \sum_{\tau =a+1}^{b}C_{l}\left( \tau
\right) \left( \nabla _{\left( a+1\right) \ast }^{\mu _{l}}f\left( \tau
\right) \right) ^{2}\right) .
\end{equation*}%
That is 
\begin{equation*}
\left( \sum_{j=a+m_{k}}^{b}\left| f\left( j\right) \right| ^{r}\right) ^{%
\frac{2}{r}}\leq \left( \sum_{j=a+m_{l}}^{b}\left| f\left( j\right) \right|
^{r}\right) ^{\frac{2}{r}}\leq 
\end{equation*}%
\begin{equation}
\delta ^{\ast }\rho ^{\ast }\left( \sum_{\tau =a+1}^{b}C_{l}\left( \tau
\right) \left( \nabla _{\left( a+1\right) \ast }^{\mu _{l}}f\left( \tau
\right) \right) ^{2}\right) =\delta ^{\ast }\rho ^{\ast }B_{l}\text{, \ for }%
l=1,...,k.  \tag{61}  \label{61}
\end{equation}%
Hence 
\begin{equation}
\left\| f\right\| _{r,[a+m_{k},b]}^{2}\leq \delta ^{\ast }\rho ^{\ast
}\left( \frac{\sum_{l=1}^{k}B_{l}}{k}\right) ,  \tag{62}  \label{62}
\end{equation}%
proving the claim.
\end{proof}


\begin{thebibliography}{9}
\bibitem{1} G. Anastassiou, \textit{Discrete fractional Calculus and
inequalities}, submitted, 2009.

\bibitem{2} D.R. Anderson, \textit{Taylor Polynomials for nabla Dynamic
equations on time scales}, Panamer. Math. J., 12(4):17-27, 2002.

\bibitem{3} F. Atici and P. Eloe, \textit{Initial value problems in discrete
fractional calculus}, Proc. AMS, 137, no. 3 (2009), 981-989.

\bibitem{4} F. Atici and P. Eloe, \textit{Discrete fractional calculus with
the nabla operator}, Electronic J. of Qualitative Theory of Differential
Equations, Spec. Ed. I, 2009, No. 1, 1-99;
http://www.math.u-szeged.hu/ejqtde/.
\end{thebibliography}
\end{document}